# Subfactors of index less than 5, part 4: vines

David Penneys
James E. Tener

URLs:
http://math.berkeley.edu/~dpenneys and http://math.berkeley.edu/~jtener
Email: dpenneys@math.berkeley.edu  jtener@math.berkeley.edu

**Abstract**  We eliminate 38 infinite families of possible principal graphs as part of the classification of subfactors up to index 5. A number-theoretic result of Calegari-Morrison-Snyder, generalizing Asaeda-Yasuda, reduces each infinite family to a finite number of cases. We provide algorithms for computing the effective constants that are required for this result, and we obtain 28 possible principal graphs. The Ostrik $d$-number test and an algebraic integer test reduce this list to 7 graphs in the index range $(4, 5)$ which actually occur as principal graphs.

## 1 Introduction

This is the fourth and final paper in the classification of subfactors to index 5. This extends the previous classifications to $3 + \sqrt{3}$ by [Jon83, Ocn88, GdlHJ89, Izu91, Pop93, Pop94, Haa94, EK98, Bis98, AH99, KO02, AY09, BMPS09]. This paper uses number theoretic techniques of [dBG91, CG94, ENO05, Ost09, CMS10] to eliminate 38 infinite families of possible principal graphs.

The following theorem is the main result of this series of papers.

**Theorem 1.1**  *There are exactly ten subfactor planar algebras other than Temperley-Lieb with index between 4 and 5: the Haagerup planar algebra and its dual [AH99], the extended Haagerup planar algebra and its dual [BMPS09], the Asaeda-Haagerup planar algebra and its dual [AH99], the 3311 Goodman-de la Harpe-Jones planar algebra and its dual [GdlHJ89, Oka91], and Izumi's self-dual 2221 planar algebra [Izu01] and its complex conjugate.*

Up to this point in this series, we have proved the following.

**Theorem 1.2**  (From [MS10, MPPS10, IJMS])  *Any exception to the above result must be a translate of one of an explicit finite list of graph pairs (which we call the vines, see Definition 2.1 and Table A).*

The point of this paper is to automate a theorem of Calegari-Morrison-Snyder [CMS10] to obtain results like those of Bisch [Bis98] and Asaeda-Yasuda [AY09] on a large scale and in a uniform manner. We apply this machinery to eliminate all 38 vines from [MS10] in the main theorem below. The relevant notation is explained at the beginning of Subsection 2.1.

**Theorem 1.3** *Suppose $(\Gamma, v)$ is a vine from the list in Table A, and $\Gamma_{|\Gamma|+j}$ is the translation of $\Gamma$ at $v$ by $j$. If $\Gamma_{|\Gamma|+j}$ is the principal graph of a non-$A_\infty$ subfactor with index in $(4, 5)$, then $\Gamma, j$ are in the following table ($v$ is always the leftmost vertex). Corresponding subfactors have been constructed for each allowed graph.*

| # | Vine $\Gamma$ | Translates | Constructed in |
|---|---|---|---|
| 4 | 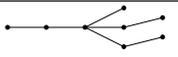 | $j = 0$ | [Izu01] |
| 7 | 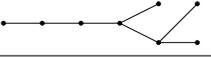 | $j = 0$ and $4$ | [AH99] and [BMPS09] |
| 9 | 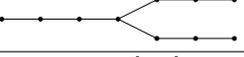 | $j = 0$ and $4$ | [AH99] and [BMPS09] |
| 21 | 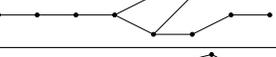 | $j = 2$ | [AH99] |
| 25 | 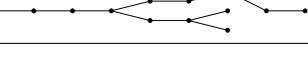 | $j = 2$ | [AH99] |

As mentioned above, the primary tool for the proof is Theorem 1.0.3 of [CMS10]. They prove that given a vine $(\Gamma, v)$, one can compute an $N(\Gamma) \in \mathbb{N}$ such that $\Gamma_n$ is not a principal graph for all $n > N(\Gamma)$, where $\Gamma_n$ is the translation of $\Gamma$ at $v$ with $n$ vertices. A similar approach was used in [AY09] to eliminate translates of the Haagerup family (Vines 7 and 9) as possible principal graphs.

The structure of this paper is as follows:

In Section 2, we give the background necessary for this paper. Subsection 2.1 explains how we obtained the list of vines in Table A. In Subsection 2.2, we recall material from [CMS10] used to calculate effective constants to eliminate the vines listed in Table A. Finally, Subsection 2.3 recalls material from [Ost09] on $d$-numbers.

In Section 3, we give algorithms for explicitly computing $N(\Gamma)$, along with algorithms for the cyclotomic test (Algorithm 3.6) and the Ostrik $d$-number test (Algorithm 3.7) used to eliminate $\Gamma_n$ for most $n \leq N(\Gamma)$. Using these algorithms, we prove Theorem 1.3 in Subsection 3.3. Necessary data for the proof is found in Tables A, B, and C.

Bundled with the arXiv source of this article are two Mathematica notebooks, named `EliminatingVinesContent.nb` and `EliminatingVinesCode.nb`, which contain all relevant calculations for what follows. These make use of a package called `FusionAtlas`; see [MS10] for a terse tutorial on its use. The most computationally intensive step in our calculations is determining if a certain equation holds for an $n$th root of unity for all $n$ less than an explicit, large bound. While each computation can be done by hand, the number of calculations necessitates the use of a computer.

We would like to thank Scott Morrison, Emily Peters, and Noah Snyder for many helpful conversations, for proofreading the manuscript, and for useful programs in the `FusionAtlas`. We would also like to thank Vaughan Jones for hosting several Bodega Bay Planar Algebra Programming Camps. Both authors would like to acknowledge support from NSF grants DMS 0401734 and DMS 0856316 and from DOD-DARPA grant HR0011-11-1-0001. The second author was also supported by an NSF Graduate Research Fellowship.

## 2 Background

### 2.1 Vines

**Definition 2.1** *A vine is a pair $(\Gamma, v)$ where $\Gamma$ is a finite, connected, bipartite graph with $|\Gamma|$ vertices and $v$ is a vertex of $\Gamma$.*

Note that [MS10] defined a vine as a pair of graphs with dual data. Our definition differs as the obstructions we use in this paper only deal with one graph at a time without dual data.

**Definition 2.2** *For a vine $(\Gamma, v)$, let $\Gamma_n$ denote the sequence of graphs obtained by adding a 2-valent tree of length $n - |\Gamma|$ to $\Gamma$ at $v$. When $\Gamma$ is a principal graph and $v$ is the initial vertex, then the initial vertex of $\Gamma_n$ is the vertex at the end of the attached 2-valent tree, i.e., $\Gamma_n$ is the translation of $\Gamma$ at $v$ by $n - |\Gamma|$. For example, $\Gamma$ translated by $j$ at $v$ is denoted $\Gamma_{|\Gamma|+j}$.*

A vine $(\Gamma, v)$ gives an infinite family $(\Gamma_n)$ of possible principal graphs of subfactors. We say a vine $(\Gamma, v)$ has been *eliminated* if we can reduce this infinite family to only finitely many possible principal graphs. When we refer to a (numbered) vine from Table A, the distinguished vertex $v$ is the leftmost vertex.

**Example 2.3** *Haagerup's classification of principal graphs of subfactors to index $3 + \sqrt{3}$ shows that the principal graph pair must be a translate of one of the following pairs of vines:*

- *Vines 7* 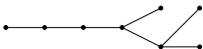 *and 9* 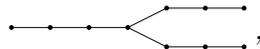 *,*
- *Vines 10* 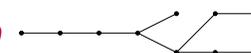 *and 11* 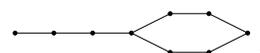 *, and*
- *Vines 21* 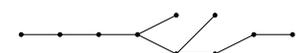 *and 25* 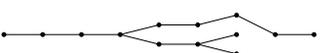 *.*

Asaeda-Yasuda eliminated Vine 9, and hence the first pair, in [AY09] using number theory.

Bisch eliminated Vine 11, and hence the second pair, in [Bis98] by showing the nonexistence of a consistent set of fusion rules. In his classification, Haagerup announced the elimination of the third vine pair. Part two of this series [MPPS10] includes a proof of this result; we also eliminate these vines in Theorem 1.3.

Each of the results mentioned in the above example used different techniques. To uniformly eliminate the vines in Table A, we use results of Calegari-Morrison-Snyder [CMS10].

### 2.2 Background from Calegari-Morrison-Snyder

By [dBG91, CG94, ENO05], the dimension of any object in a fusion category must be an algebraic integer in a cyclotomic field. The even part of the standard invariant of a finite-depth subfactor $A \subset B$ is a fusion category; hence the index $[B:A] = \dim(_A L^2(B)_A)$ is a cyclotomic integer. By [Jon86], if $A \subset B$ is finite-depth with principal graph $\Gamma$, then $[B:A] = \|\Gamma\|^2$, and thus the square of the Frobenius-Perron eigenvalue of $\Gamma$ is a cyclotomic integer. In fact, more is true; the square of every multiplicity-one eigenvalue of the adjacency matrix of $\Gamma$ is a cyclotomic integer by Lemma 3.0.7 of [CMS10]. Using this result, combined with some number-theoretic techniques [Cas69, Lox72, GHM09], Calegari-Morrison-Snyder obtained the following theorem.

**Theorem 2.4** (From [CMS10]) *Given a vine $(\Gamma, v)$ such that $\Gamma_n$ is not $A_n$ or $D_n$ for $n > |\Gamma|$, there are constants $K(\Gamma)$ and $|R|$ which can be effectively computed directly from $(\Gamma, v)$ such that $\Gamma_n$ is not a principal graph whenever*

$$n > 4K(\Gamma) + 9|R|. \tag{2.1}$$

This paper provides the machinery for implementing the above theorem on the large scale required for recent classification results. Table A contains a list of 38 vines and the constants required for the application of Theorem 2.4. Table A is obtained from the list $\mathcal{V}_\infty$ in Theorem 6.1 of [MS10]. The list $\mathcal{V}_\infty$ contains pairs of graphs with dual data; we forget the dual data, uncouple the graphs, remove duplicates and translates, and order the vines by increasing depth.

We use Mathematica to automate the algorithms given in Section 3 to calculate $K(\Gamma)$ and an upper bound $\mathcal{R}$ on $|R|$ directly from $\Gamma$. Hence $\Gamma_n$ is not a principal graph whenever $n$ is larger than

$$N(\Gamma) = 4K(\Gamma) + 9\mathcal{R} \geq 4K(\Gamma) + 9|R|. \tag{2.2}$$

In practice, we find $\|\Gamma_n\|^2$ is not cyclotomic for most $n$ much smaller than $N(\Gamma)$. This leads to the interesting question of whether the bounds in Equations (2.1) and (2.2) can be improved.

We recall material from [CMS10] which is essential in calculating the constants $K(\Gamma)$ and $\mathcal{R}$ from Equation (2.2). The notation used here will be used throughout Section 3.

Let $(\Gamma, v)$ be a vine such that $\Gamma_n$ is not $A_n$ or $D_n$ for $n > |\Gamma|$. Let $M_n$ denote the adjacency matrix of $\Gamma_n$, and let $P_n \in \mathbb{Z}[x]$ denote the characteristic polynomial of $M_n$. Note $\deg(P_n) = n$ as $\Gamma_n$ has $n$ vertices.

**Remark 2.5** *As $\Gamma_n$ is bipartite, $r \in \mathbb{R} \setminus \{0\}$ is an eigenvalue of $M_n$ if and only if $-r$ is, in which case they occur with the same multiplicity.*

**Lemma 2.6** (From [CMS10]) *Let $x = t + t^{-1}$, and write $P_n(x) = F_n(t) \in \mathbb{Z}[t, t^{-1}]$.*

(1) *The matrix $M_n$ is symmetric and the roots of $P_n(x)$ are all real.*

(2) *The polynomials $P_n$ satisfy the recurrence*

$$P_n(x) = xP_{n-1}(x) - P_{n-2}(x). \tag{2.3}$$

(3) *There is a fixed Laurent polynomial $A \in \mathbb{Z}[t, t^{-1}]$ such that*

$$F_n(t)\left(t - t^{-1}\right) = t^n A(t) - t^{-n} A(t^{-1}). \tag{2.4}$$

**Remark 2.7** *If $A \in \mathbb{Z}[t, t^{-1}]$, then the sequence $[t^n A(t) - t^{-n} A(t^{-1})]_{n \in \mathbb{N}}$ satisfies Recurrence (2.3).*

We use the letter $\lambda$ to refer to a root of $P_n(x)$ and the letter $\rho$ to refer to the corresponding roots of $F_n(t)$, where $\lambda = \rho + \rho^{-1}$.

**Lemma 2.8** (From [CMS10]) *Let $K(\Gamma) = \sum \rho^4$ such that $\rho$ is a root of $F_n$ for $n$ sufficiently large compared to $\deg(A)$. The constant $K(\Gamma)$ is well defined since the sum of the 4th powers of the roots of $F_n$ depends only on the first four coefficients of $F_n$, which is independent of $n$ for $n$ sufficiently large.*

Algorithm 3.2 computes $K(\Gamma)$ and determines which $n$ are "sufficiently large" in the preceding lemma. We now define $R$.

**Definition 2.9** (From [CMS10]) *Let $R = R_1 \cup R_2 \cup R_3$, where the $R_i$'s are the (multi)sets of roots <u>with multiplicity</u> of $P_n(x)$ given by:*

(1) *$R_1$ is the roots of the form $\zeta + \zeta^{-1}$, where $\zeta$ is a root of unity.*

(2) *$R_2$ is the set of roots which appear with multiplicity $\geq 2$.*

(3) *$R_3$ is the set of roots equal to $\lambda^2 - 2$ for $\lambda = 1 + 2\cos(2\pi/7)$ or $2\cos(\pi/30) + 2\cos(13\pi/30)$.*

Algorithm 3.4 computes $\mathcal{R}$, an upper bound for $|R|$ independent of $n$. The following definitions and results will be useful for bounding the size of $R_2$.

**Definition 2.10** *A polynomial $B \in \mathbb{Z}[t]$ with nonzero constant term is called <u>self-reciprocal</u> if $B(t) = \pm t^{\deg(B)} B(t^{-1})$.*

*Remark.* The minimal polynomial over $\mathbb{R}$ for $r \in \mathbb{C}$ is self-reciprocal when $|r| = 1$.

**Lemma 2.11** (From [CMS10]) *Factor $A(t) = t^{-s} B(t) C(t)$, where $s \geq 0$, and $B, C$ are polynomials with nonzero constant term such that $B$ is a maximal self-reciprocal polynomial factor of $A$ (so $C$ has no roots on the unit circle). Then there is a $d > s$ such that whenever $n \geq d$,*
$$(2(n-s) + \deg(B))|C(t)| - |C'(t)| - |C'(t^{-1})| > 0, \qquad (2.5)$$
*so every repeated root of $F_n$ on the unit circle is a root of $B$.*

In practice, we calculate this $d$ by increasing $n$ until Inequality (2.5) is satisfied, which we verify numerically. This is only necessary for some vines (see Subsection 3.2), and we compute an upper bound for $d$ for these vines in Table C.

Finally, we will make use of:

**Theorem 2.12** (Descartes' rule of signs) *Suppose $A \in \mathbb{R}[t, t^{-1}]$. List the coefficients of $A$ by descending power of $t$ (excluding zeroes), and let $\mathrm{SignChanges}(A)$ be the number of sign changes in the list. Let $r$ be the number of positive roots of $A$, counted with multiplicity. Then $r \leq \mathrm{SignChanges}(A)$ and $r = \mathrm{SignChanges}(A) \mod 2$.*

## 2.3 Background on $d$-numbers

The algorithms of Section 3 are not sufficient to eliminate all of the graphs obtained from the vines in Table A. We recall material from [Ost09] to eliminate some exceptional graphs in Subsection 3.3.

**Definition 2.13** *An algebraic integer $\alpha$ is called a <u>$d$-number</u> if the following condition holds:*

*Let $p(x) = x_n + a_1 x_{n-1} + \cdots + a_n$ be the minimal polynomial of $\alpha$ over $\mathbb{Q}$ (so $a_i \in \mathbb{Z}$). Then for any $i = 1, \ldots, n$, the number $(a_i)^n$ is divisible by $(a_n)^i$.*

**Definition 2.14** *Given a fusion category $\mathcal{C}$, its global dimension is*

$$\sum_{\text{simple } X \in \mathcal{C}} \dim_{FP}(X)^2$$

*where $\dim_{FP}$ is the Frobenius-Perron dimension.*

**Theorem 2.15** (Corollary 1.3 of [Ost09]) *The global dimension of a fusion category is a $d$-number.*

Algorithm 3.7, the Ostrik $d$-number test, uses this theorem as an obstruction to possible principal graphs.

## 3 Algorithms

### 3.1 The general case

In this section, we give algorithms for computing $K(\Gamma)$ and $\mathcal{R}$ for vines. These algorithms rely on basic linear algebra and results from [CMS10, GHM09]. We also explain the cyclotomic test and Ostrik $d$-number test, two powerful obstructions for possible principal graphs.

We will use the notation and results of Subsection 2.2. All algorithms will be applied to a fixed vine $(\Gamma, v)$.

**Definition 3.1** *By Equation (2.4), if $A$ is written in descending monomial order*

$$A(t) = a_r t^r + \cdots + a_{-s} t^{-s},$$

*and if $n > s \geq 0$, then the analytic part of $(t - t^{-1}) F_n(t)$ is given by $t^n A(t)$ and the principal part of $(t - t^{-1}) F_n(t)$ is given by $-t^{-n} A(t^{-1})$. For $n > s$, we say $(t - t^{-1}) F_n(t)$ is separated, and $s$ is the separation constant for $\Gamma$.*

This separation allows us to calculate $K(\Gamma)$, as if $n \geq s$ then $n$ is "sufficiently large" for Lemma 2.8. We use this to give an algorithm to calculate $A$, $K(\Gamma)$ and the separation constant $s$.

**Algorithm 3.2** (Separation)

(1) We begin by computing $A$:

  (a) Set $k = |\Gamma| + 1$.
  (b) Let $C_k$ be the Laurent polynomial obtained from the analytic part of $G_k$ by multiplying by $t^{-k}$. We call $C_k$ the candidate for $A$.
  (c) If $C_k \neq C_{k+1}$ (so $G_k$ is not yet separated), increase $k$ by 1 and return to (b).
  (d) Now $C_k = C_{k+1}$. If $F_{|\Gamma|+j}(t) \neq t^{|\Gamma|+j} C_k(t) - t^{-|\Gamma|-j} C_k(t^{-1})$ for $j = 1$ or $2$, increase $k$ by 1 and return to (b).
  (e) Now $C_k = C_{k+1}$ and $F_{|\Gamma|+j}(t) = t^{|\Gamma|+j} C_k(t) - t^{-|\Gamma|-j} C_k(t^{-1})$ for $j = 1$ and $2$, so $A = C_k$ by Recurrence (2.3) and Remark 2.7.

This process terminates since some $A$ exists by Lemma 2.6.

(2) Calculate $s$ as in Definition 3.1. Observe that this is the same $s$ as in Lemma 2.11.

(3) Numerically calculate $K(\Gamma) = \sum \rho^4$ such that $\rho$ is a root of $F_n$ for $n = s+1$. Note that a numerical calculation suffices as $K(\Gamma)$ is an integer.

The algorithm to compute $\mathcal{R}$ uses the following lemma, which is based on a lemma from [CMS10].

**Lemma 3.3** *Suppose $C(t) \in \mathbb{Z}[t]$ is as in Lemma 2.11 ($C$ has no self-reciprocal factors), $\gamma \in \mathbb{N}$, $c = \deg(C)$, and $\zeta_k = e^{\pi i/k}$. Set $H_n(t) := t^{2n-\gamma}C(t) \mp t^c C(t^{-1})$.*

(1) *If $H_n(\zeta_k) = 0$, then $k \leq 2Lc$ where $L$ is the product of primes $p$ less than or equal to $2(\#\text{monomial summands of } C(t))$. Note $L$ is independent of $n$.*

(2) *$H_n(\zeta_k) = 0$ for some $n$ if and only if $H_m(\zeta_k) = 0$ for some $1 \leq m \leq \hat{k}$ where $\hat{k} = k$ if $k$ is odd and $\hat{k} = k/2$ if $k$ is even.*

(3) *Let $S = \{k : \zeta_k \text{ is a root of some } H_n\}$. The number of roots of unity $\zeta$ satisfying $H_n(\zeta) = 0$ is periodic in $n$ with period $\text{LCM}\{\hat{k} \mid k \in S\}$.*

**Proof**

(1) This follows immediately from Theorem 2.1 of [GHM09].

(2) Suppose that $H_n(\zeta_k) = 0$. Then $H_m(\zeta_k) = 0$ if and only if $k$ divides $2(n-m)$ if and only if $n = m \mod \hat{k}$.

(3) Note that $S$ is finite by (1). The rest follows from (2). □

**Algorithm 3.4** *Calculate $\mathcal{R}$ as follows:*

(1) Factor $A(t) = t^{-s}B(t)C(t)$ as in Lemma 2.11.

(2) Compute an upper bound on $d$ from Lemma 2.11 by finding an $n$ for which Inequality (2.5)
$$(2(n-s) + \deg(B))|C(t)| - |C'(t)| - |C'(t^{-1})| > 0$$
holds. This requires human intervention: we numerically graph the function on the left hand side of the inequality (see Section 3.2 and Table C).

(3) Set $\alpha = n - s - b + c$ and $\beta = 2(n-s) + b + c$ so that
$$(t - t^{-1})F_n(t) = t^\alpha B(t)\left(t^\beta C(t) \mp t^c C(t^{-1})\right).$$

(4) Solve $B(t) = 0$.

  (a) Let $r_1$ be the number of complex conjugate pairs of solutions which are roots of unity, as complex conjugates yield the same root of $P_n(x)$.

  (b) Let $r_2$ be the number complex conjugate pairs of solutions which are repeated roots.

(5) We now bound the number of roots of unity which are zeroes of
$$H_n(t) := t^\beta C(t) \mp t^c C(t^{-1})$$
for some $n \geq s$, but not equal to $\pm 1$. We use Lemma 3.3 with $\gamma = b + c - 2s$.

   (a) Compute the $L$ of Lemma 3.3.
   (b) Compute $S = \{k \mid \zeta_k \text{ is a root of some } H_n\}$ as follows: for all $k \leq 2Lc$, check if $H_n(\zeta_k) = 0$ for $3 \leq n \leq \hat{k}$. In fact, by Remark 10.1.8 of [CMS10], we need only check $k$ such that $k$ divides $mL$ for some $m \leq 4c$.
   (c) Let $\ell = \text{LCM}\{\hat{k} \mid k \in S\}$.
   (d) For $i = 1, \ldots, \ell$, let $r_{3,i}$ be the number of roots of unity which are roots of $F_i$. Set $r_3 = \max\{r_{3,i} \mid i = 1, \ldots, \ell\}$.

(6) Set $r_4 = 2\,\text{SignChanges}(A) + 1$

(7) Set $\mathcal{R} = r_1 + r_2 + r_3 + r_4$.

**Lemma 3.5** *The number $\mathcal{R}$ calculated in Algorithm 3.4 satisfies $\mathcal{R} \geq |R|$ when $n \geq d > s$.*

**Proof** First, by Lemma 3.3, note that $r_1 + r_3 \geq |R_1 \setminus \{\pm 1\}|$, where we treat these sets with multiplicity (we remove all occurrences of $\pm 1$ from $R_1$).

As $n \geq d$, Lemma 2.11 assures us that all repeated roots of $F_n(t)$ on the circle are in fact roots of $B(t)$, so $r_2$ bounds the number of repeated roots of $P_n(x)$ corresponding to roots of $F_n(t)$ on the circle.

Suppose now that $x_0 = t_0 + t_0^{-1}$ is a repeated root of $P_n(x)$ for $t_0 \in \mathbb{R}$. Then

   (1) $-x_0$ is a root of $P_n(x)$ with multiplicity $m$ by Remark 2.5, and
   (2) $t_0$, $t_0^{-1}$ and $-t_0$, $-t_0^{-1}$ are all roots of $F_n(t)$ with multiplicity $m$ by the definition of $F_n(t)$.

Hence to count all repeated roots of $P_n(x)$ which are of the form $x_0 = t_0 + t_0^{-1}$ for $t_0 \in \mathbb{R}$, it suffices to count only the repeated roots of $F_n(t)$ in $[-1, 1]$. This is equal to the number of positive repeated roots of $F_n(t)$ plus the multiplicity of $-1$ as a root of $F_n(t)$.

We use Descartes' rule of signs to overcount the positive repeated roots of $F_n(t)$ along with the exceptional roots $R_3$. From Definition 3.1, for all $n \geq s$
$$\text{SignChanges}((t - t^{-1})F_n(t)) = 2\,\text{SignChanges}(A(t)) + 1,$$
so $F_n(t)$ has at most $2\,\text{SignChanges}(A) = r_4 - 1$ repeated positive roots on the real line.

Now there are 3 cases depending on the multiplicity of $-1$ as a root of $F_n(t)$, which is equal to the multiplicity of $-2$ as a root of $P_n(x)$. However, recall $\pm 2$ occur with the same multiplicity for $P_n(x)$ by Remark 2.5, so $\pm 1$ occur with the same multiplicity for $F_n(t)$.

<u>Case 1:</u> Suppose $\pm 1$ are not roots of $F_n(t)$. Hence $r_1 + r_3 \geq |R_1|$, and $r_2 + r_4 \geq |R_2 \cup R_3|$.

<u>Case 2:</u> Suppose $\pm 1$ are roots of $F_n(t)$ with multiplicity $> 1$. Then $\pm 1$ are roots of $B$ by Lemma 2.11. Hence $r_1 + r_2 + r_3 + r_4 \geq |R_1 \cup R_2 \cup R_3|$.

**Case 3:** Suppose $\pm 1$ are roots of $F_n(t)$ with multiplicity 1. Then $\pm 1$ is not necessarily a root of $B$. Recall that $r_4 = 2\operatorname{SignChanges}(A) + 1$ overcounts the positive real roots of $F_n$ by 1, which accounts for the possibility that $-1$ is not a root of $B$. So once again, $r_1 + r_2 + r_3 + r_4 \geq |R_1 \cup R_2 \cup R_3|$.

Hence $\mathcal{R} = r_1 + r_2 + r_3 + r_4 \geq |R_1 \cup R_2 \cup R_3| = |R|$. □

*Remark.* Since $\lambda_i^2 - 2 > 3$ for $i = 1, 2$, any possible subfactor with principal graph that could have $\lambda_i^2 - 2$ as a root of $P_n(x)$ must have index greater than 9. Hence $R_3 = \emptyset$ for all vines considered in this paper.

*Remark.* Note that if $A$ has only one sign change, then we may skip Steps 2, 4b, and 6 in Algorithm 3.4 (see Subsection 3.2), and $\mathcal{R} = r_1 + r_3$.

After computing $K(\Gamma)$ and $\mathcal{R}$, we set $N(\Gamma) = 4K(\Gamma) + \mathcal{R}$ by Equation (2.2). To eliminate most of the $\Gamma_n$ for $n \leq N(\Gamma)$, we use the following:

**Algorithm 3.6** (Cyclotomic Test) *Let $\alpha$ be an algebraic integer and let $M > 0$.*

(1) *Calculate the minimal polynomial $P(x)$ of $\alpha$.*

(2) *Calculate the list $L$ of the primes $\leq M$ which do not divide the discriminant of $P(x)$.*

(3) *If $L$ is empty,*

    (a) *The algorithm terminates, and $\alpha$ passes the cyclotomic test with upper bound $M$.*

    (b) *Otherwise, pick the smallest $p \in L$. If $P(x)$ does not have uniform degree irreducible factors mod $p$,*

        i. *Then $\alpha$ is not a cyclotomic integer by Theorem 4.6 of [Nar04], and we say $\alpha$ fails the cyclotomic test for prime $p$.*

        ii. *Otherwise, replace $L$ with $L \setminus \{p\}$ and return to (3).*

*Remark.* If $\Gamma$ is a bipartite graph and $\|\Gamma\|^2$ fails the above test, then $\Gamma$ is not a principal graph by [ENO05, Corollary 8.54]. Some exceptional graphs which are not principal graphs pass the cyclotomic test with $M = 200$ and presumably have cyclotomic square norm (see Tables A and B). Hence the cyclotomic test is not sufficient for the proof of Theorem 1.3, and we need a few alternate obstructions.

**Algorithm 3.7** (Ostrik $d$-number Test) *Suppose $(\Gamma, *)$ is a bipartite graph with distinguished even vertex $*$.*

(1) *Calculate the Frobenius-Perron dimensions of $\Gamma$.*

(2) *Calculate $\sum_{\text{even } v} \dim(v)^2$, the global even dimension of $\Gamma$.*

(3) *If it is not a $d$-number, then $(\Gamma, *)$ fails the Ostrik $d$-number test.*

*Remark.* If $(\Gamma, *)$, a bipartite graph with distinguished even vertex, fails the above test, then $(\Gamma, *)$ is not a principal graph by Theorem 2.15. There are examples of graphs which are not principal graphs which pass both the cyclotomic test for $M = 200$ and the Ostrik $d$-number test, e.g., Vine 32 translated by 1 (see Remark 3.14).

## 3.2 The Salem case

Algorithm 3.4 will compute $\mathcal{R}$ for any vine. In practice, the bound is often bad, and a modest simplification produces a better result with fewer and faster computations. Consider the following:

**Example 3.8** Let $\Gamma$ be Vine 18: 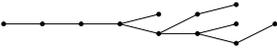. We calculate that

$$A(t) = -\frac{1}{t^{15}} - \frac{2}{t^{13}} - \frac{4}{t^{11}} - \frac{6}{t^9} - \frac{7}{t^7} - \frac{6}{t^5} - \frac{3}{t^3} + t$$
$$= t^{-15}\left(t^2+1\right)\left(t^{14} - t^{12} - 2t^{10} - 4t^8 - 3t^6 - 3t^4 - t^2 - 1\right).$$

Hence $C(t) = t^{14} - t^{12} - 2t^{10} - 4t^8 - 3t^6 - 3t^4 - t^2 - 1$, and one calculates (numerically) that the $d$ of Lemma 2.11 is greater than 245.

This example motivates partitioning the vines based on Salem numbers.

**Definition 3.9** A real algebraic integer $\alpha > 1$ is called a <u>Salem number</u> if all its Galois conjugates have modulus $\leq 1$, and at least one conjugate has modulus $1$.

If $A(t) = t - D(t^{-1})$ where $D \in \mathbb{Z}[t]$ is a polynomial with nonnegative coefficients, then for $n > s$, the separation constant, $F_n(t)$ has exactly two sign changes. By Descartes' rule of signs 2.12, $F_n(t)$ has either 0 or 2 positive roots. We know that $P_n(x)$ has a unique positive root $x_n = t_n + t_n^{-1} > 2$ (the Froebenius-Perron eigenvalue), so $t_n, t_n^{-1}$ are the two positive roots of $F_n(t)$. Now $P_n(x)$ has only real roots, so the other roots of $F_n(t)$ lie on the circle. Without loss of generality, $t_n > 1 > t_n^{-1}$, and $t_n$ is often Salem, with some small exceptions. Now all Galois conjugates of $t_n$ are roots of multiplicity one, and all other roots of $F_n(t)$ are roots of unity by the following.

**Fact 3.10** Suppose $|r| = 1$, and all of $r$'s Galois conjugates lie on the circle. Then $r$ is a root of unity.

Thus as long as $n > s$ (and $\|\Gamma_n\|$ is not an exceptional root in $R_3$), $R = R_1$ as all repeated roots must be roots of unity. This discussion motivates the following definition.

**Definition 3.11** A vine is called <u>Salem</u> if $A(t) = t - D(t^{-1})$ where $D \in \mathbb{Z}[t]$ has nonnegative coefficients.

Of the vines in Table A, only 8 are <u>not</u> Salem: 22, 26, 27, 29, 33, 34, 36, and 37.

For the non-Salem vines, a simple numerical calculation in Table C shows that Inequality (2.5) is satisfied if $n \geq 200$. Hence $200 \geq d$ from Lemma 2.11, and since $N(\Gamma) > 200$ for all non-Salem vines, we do not include $d$ in Table A. (In fact, we calculated $70 \geq d$, but the graphs in Table C are much clearer using a worse upper bound.)

## 3.3 Proof of Theorem 1.3

**Theorem 3.12** *Suppose $\Gamma$ is one of the 38 vines in the second column of Table A. Then columns 3 and 4 give the $s$ and $K(\Gamma)$ computed by Algorithm 3.2, column 5 gives $\mathcal{R}$ computed by Algorithm 3.4, and column 6 gives $N(\Gamma) = 4K(\Gamma) + 9\mathcal{R}$ obtained from Equation (2.2). Thus $\Gamma_n$ is not a principal graph for all $n > N(\Gamma)$.*

*For all $j + |\Gamma| \leq N(\Gamma)$, except for those $j$'s appearing in the final "Exceptions" column, $\Gamma_{|\Gamma|+j}$ fails the cyclotomic test for one of the primes $p$ given in column 6.*

*Hence starting with the 38 vines appearing in Table A, only 28 exceptional graphs can be principal graphs.*

**Proposition 3.13** *Of the 28 exceptional graphs, the 15 graphs in Table B fail the Ostrik $d$-number test (Algorithm 3.7). The minimal polynomials of the global even dimensions of the graphs are given in column 3.*

**Remark 3.14** *The exceptional graph given by Vine 32 translated by $j = 1$ is not a principal graph. Morrison and Ostrik have shown that for any fusion ring consistent with the even part of this graph, there is a formal codegree which does not lie in the field generated by the Frobenius-Perron dimensions (this cannot happen by a corollary of Lemma 3.1 in [Ost09]). This proof will appear in a future paper. Note that the norm squared of Vine 32 translated by $j = 1$ is exactly 5, so this argument is not essential to the proof of Theorem 1.3.*

**Proposition 3.15** *Consider Vines 1, 2, and 3 translated by $j = 0$. For these three graphs, the dimension of $V_{3,1}$, the bottom vertex at depth 3, is not an algebraic integer. Thus these graphs are not principal graphs of subfactors (all the Frobenius-Perron dimensions of a multi-fusion ring are graph norms, and hence algebraic integers). The pertinent information is listed in the following table:*

| Vine # | Minimal polynomial of $\dim(V_{3,1})$ |
|:---:|:---:|
| 1 | $3x^4 - 8x^2 + 1$ |
| 2 | $2x^4 - 18x^2 + 3$ |
| 3 | $2x^4 - 18x^2 + 3$ |

**Proof of Theorem 1.3** From Theorem 3.12, Propositions 3.13 and 3.15, and Remark 3.14, we have exactly 9 graphs which remain as possible principal graphs, seven of which are listed in Theorem 1.3. The remaining two are Vines 5 and 6 translated by one:

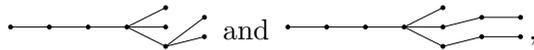 and ,

which are principal graphs of subfactors at index 5 corresponding the the inclusion of groups $A_4 \subset A_5$. □

# A  Table of vines and effective constants

| # | Vine | s | $K(\Gamma)$ | $\mathcal{R}$ | $N(\Gamma)$ | $p$ | Exceptions |
|---|------|---|---|---|---|---|---|
| 1 | 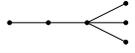 | 6 | 10 | 4 | 76 | 179,181 | $j=0,1$ |
| 2 | 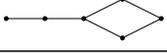 | 6 | 12 | 3 | 75 | 41,43 | $j=0,1$ |
| 3 | 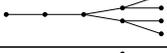 | 12 | 10 | 6 | 94 | 41,43 | $j=0,1$ |
| 4 | 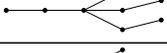 | 10 | 10 | 6 | 94 | 29,31 | $j=0$ |
| 5 | 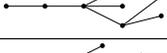 | 10 | 14 | 7 | 119 | 59,61 | $j=0$ |
| 6 | 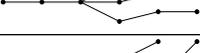 | 14 | 10 | 9 | 121 | 59,61 | $j=0$ |
| 7 | 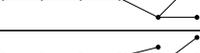 | 8 | 6 | 7 | 87 | 41,43 | $j=0,4$ |
| 8 | 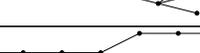 | 10 | 14 | 5 | 101 | 37,41 | None |
| 9 | 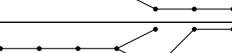 | 12 | 2 | 9 | 89 | 41,43 | $j=0,4$ |
| 10 | 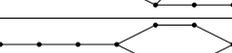 | 12 | 6 | 7 | 87 | 59 | $j=0,2$ |
| 11 | 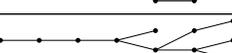 | 10 | 4 | 6 | 70 | 59 | $j=0,2$ |
| 12 | 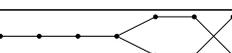 | 14 | 10 | 11 | 139 | 73,79 | $j=0$ |
| 13 | 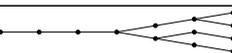 | 12 | 8 | 10 | 122 | 73,79 | $j=0$ |
| 14 | 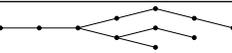 | 18 | 10 | 9 | 121 | 37,41 | None |
| 15 | 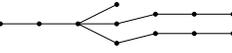 | 16 | 6 | 11 | 123 | 73,79 | $j=0,2$ |
| 16 | 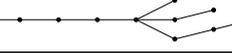 | 18 | 10 | 11 | 139 | 103,107 | None |
| 17 | 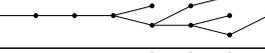 | 12 | 10 | 7 | 103 | 43,47 | None |
| 18 | 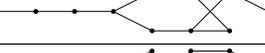 | 16 | 10 | 12 | 148 | 59,61 | None |
| 19 | 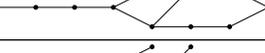 | 14 | 8 | 11 | 131 | 59,61 | None |
| 20 | 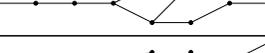 | 14 | 8 | 8 | 104 | 5,7 | None |
| 21 | 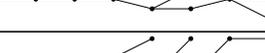 | 12 | 6 | 8 | 96 | 37,41 | $j=2$ |
| 22 | 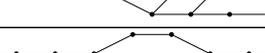 | 16 | 18 | 17 | 225 | 47,53 | None |
| 23 | 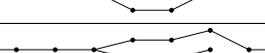 | 16 | 10 | 9 | 121 | 5,7 | None |
| 24 | 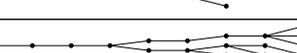 | 14 | 8 | 8 | 104 | 5,7 | None |
| 25 | 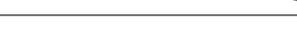 | 18 | 6 | 11 | 123 | 37,41 | $j=2$ |
| 26 | 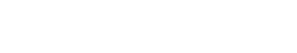 | 30 | 22 | 24 | 304 | 47,53 | None |

| 27 | | 24 | 18 | 22 | 270 | 47,53 | None |
| 28 | | 18 | 8 | 11 | 131 | 59,61 | $j=2$ |
| 29 | | 28 | 18 | 29 | 333 | 59,61 | $j=2$ |
| 30 | | 24 | 14 | 14 | 182 | 59,61 | $j=2$ |
| 31 | | 22 | 12 | 13 | 165 | 59,61 | $j=2$ |
| 32 | | 16 | 10 | 12 | 148 | 37,41 | $j=1$ |
| 33 | | 26 | 16 | 23 | 271 | 17,19 | None |
| 34 | | 30 | 20 | 27 | 323 | 17,19 | None |
| 35 | | 24 | 10 | 18 | 202 | 43,47 | None |
| 36 | | 26 | 16 | 23 | 271 | 17,19 | None |
| 37 | | 22 | 12 | 17 | 201 | 17,19 | None |
| 38 | | 28 | 6 | 20 | 204 | 43,47 | None |

# B  Table of exceptional graphs which fail the Ostrik $d$-number test

| Vine # | Translation | Minimal polynomial of global even dimension |
|:---:|:---:|:---:|
| 1 | $j = 1$ | $x^2 - 32x + 56$ |
| 2 | $j = 1$ | $x^2 - 63x + 105$ |
| 3 | $j = 1$ | $x^2 - 63x + 105$ |
| 10 | $j = 0$ | $x^2 - 65x + 275$ |
| 10 | $j = 2$ | $x^3 - 338x^2 + 2535x - 4225$ |
| 11 | $j = 0$ | $x^2 - 65x + 275$ |
| 11 | $j = 2$ | $x^3 - 338x^2 + 2535x - 4225$ |
| 12 | $j = 0$ | $x^3 - 108x^2 + 1377x - 4617$ |
| 13 | $j = 0$ | $x^3 - 108x^2 + 1377x - 4617$ |
| 15 | $j = 0$ | $5x^3 - 143x^2 + 676x - 845$ |
| 15 | $j = 2$ | $x^2 - 156x + 792$ |
| 28 | $j = 2$ | $x^3 - 684x^2 + 8505x - 26163$ |
| 29 | $j = 2$ | $x^3 - 684x^2 + 8505x - 26163$ |
| 30 | $j = 2$ | $x^3 - 684x^2 + 8505x - 26163$ |
| 31 | $j = 2$ | $x^3 - 684x^2 + 8505x - 26163$ |

# C  Table of numerical calculations for $d$'s of non-Salem vines

The second column of the following table is the graph of the left hand side of (2.5):
$$(2(200-s) + \deg(B))|C(t)| - |C'(t)| - |C'(t^{-1})|$$
for $n = 200$ and $t = e^{i\theta}$ for $0 \leq \theta \leq 2\pi$ on a logarithmic scale. As the graphs are clearly positive, we have $d \leq 200$ as in Lemma 2.11. Note that $200 < N(\Gamma)$ for all non-Salem vines in Table A.

| # | Graph | # | Graph |
|---|---|---|---|
| 22 | | 33 | |
| 26 | | 34 | |
| 27 | | 36 | |
| 29 | | 37 | |

This paper is available online at `arXiv:1010.3797`.